\documentclass[Afour,sageh,times]{sagej}

\usepackage{moreverb,url}

\usepackage[colorlinks,bookmarksopen,bookmarksnumbered,citecolor=red,urlcolor=red]{hyperref}
\newcommand\BibTeX{{\rmfamily B\kern-.05em \textsc{i\kern-.025em b}\kern-.08em
T\kern-.1667em\lower.7ex\hbox{E}\kern-.125emX}}

\def\volumeyear{2016}

\usepackage{amsmath,amsfonts,amssymb}
\usepackage{amsthm}
\newtheorem{definition}{Definition}
\usepackage{graphicx}
\usepackage{braket}
\usepackage{breqn}
\usepackage{booktabs}
\usepackage{dcolumn}
\usepackage{enumitem,comment}
\usepackage[table]{xcolor}
\usepackage{blkarray}
\usepackage{bm}
\usepackage{float,titlesec}

\setcitestyle{citesep={,}}

\usepackage{algorithm}
\usepackage[noend]{algpseudocode}
\usepackage{caption}
\usepackage{subcaption}
\allowdisplaybreaks 

\definecolor{firstColor}{cmyk}{0.053333333,0.0088888888,0.0088888888,0.017777777}
\definecolor{secondColor}{cmyk}{0.01,0.01,0.11,0.02}

\begin{document}

\runninghead{Dasari, Im, and Geerhart}

\title{Complexity and mission computability of adaptive  computing systems}

\author{Venkat R. Dasari\affilnum{1} Mee Seong Im\affilnum{2} Billy Geerhart\affilnum{1}}

\affiliation{\affilnum{1}U.S. Army Research Laboratory, Aberdeen Proving Ground, MD 21005\\
\affilnum{2}U.S. Military Academy, West Point, NY 10996}

\corrauth{Venkat R. Dasari, 
Army Research Laboratory}

\email{venakteswara.r.dasari.civ@mail.mil}

%\affil[a]{U.S. Army Research Laboratory, Aberdeen Proving Ground, MD 21005}
%\affil[b]{U.S. Military Academy, West Point, NY 10996}

%\authorinfo{Further author information: (Send correspondence to Venkat R. Dasari)  \\
%Venkat R.~Dasari: E-mail: venkateswara.r.dasari.civ@mail.mil, Telephone: +01-410-278-2846  \\  
% Mee Seong Im: E-mail: meeseong.im@usma.edu, Telephone: +01-845-938-5649 %  \\ 
%}

%\email{alistair.smith@sunrise-setting.co.uk}

\begin{abstract}
There is a subset of computational problems that are computable in polynomial time for which an existing algorithm may not complete due to a lack of high performance technology on a mission field. We define a subclass of deterministic polynomial time complexity class called mission class, as many polynomial problems are  not computable in mission time. 
By focusing on such subclass of languages in the context for successful military applications,  
we also discuss their computational and communicational constraints. We investigate feasible (non)linear models that will minimize energy and maximize memory, efficiency, and computational power, and also provide an approximate solution obtained within a pre-determined length of computation time using limited resources so that an optimal solution to a language could be determined.  
\end{abstract}


\maketitle

\section{Introduction}
\label{sec:intro}  % \label{} allows reference to this section
\begin{figure}
\centering
\begin{subfigure}[b]{1.0\columnwidth} 
  %\centering
  \includegraphics[width=1.0\linewidth]{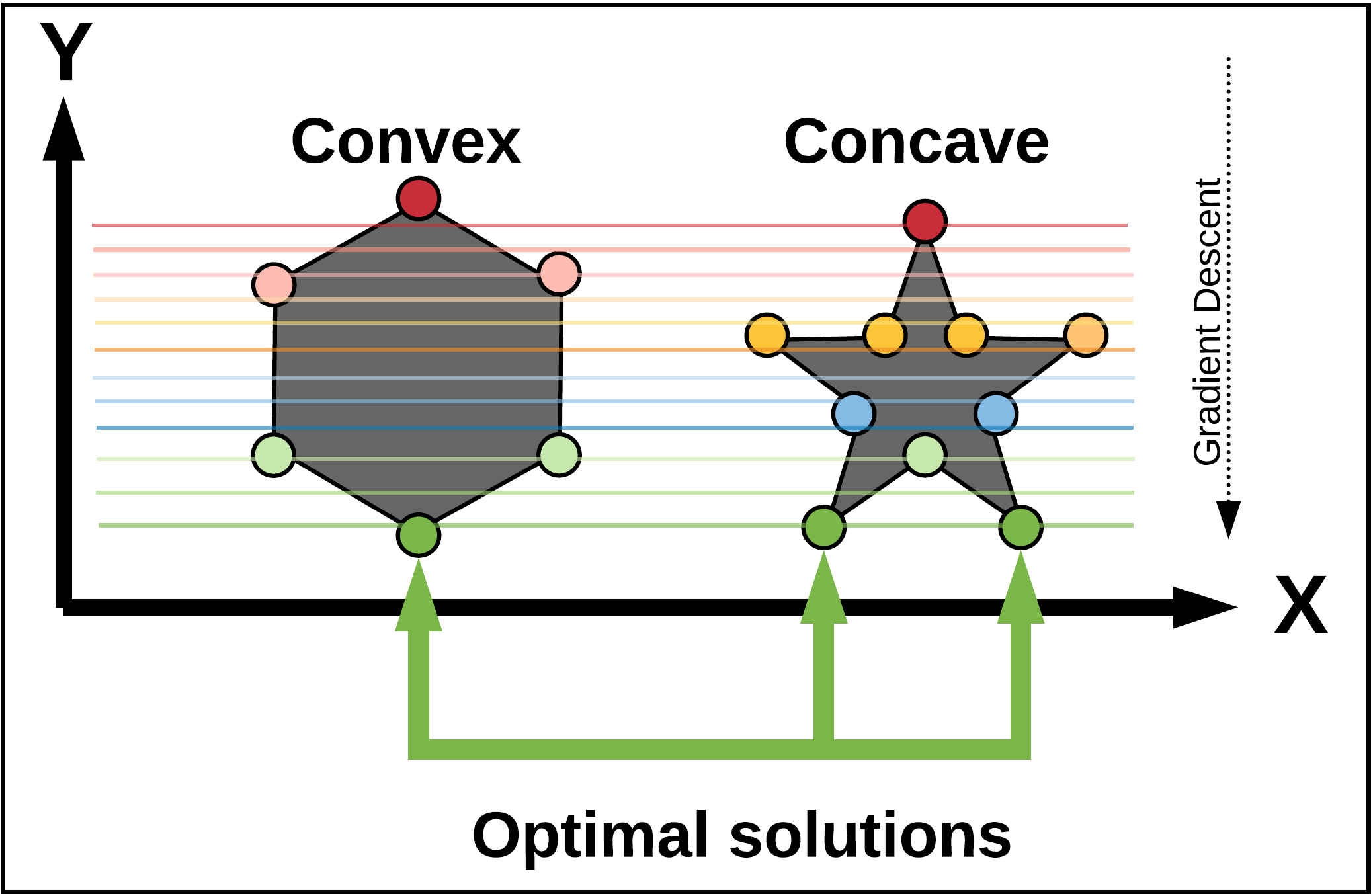}
  \caption{Convex and concave linear optimization}
  \label{fig:sub1}
\end{subfigure}%

\begin{subfigure}[b]{1.0\columnwidth}
  %\centering
  \includegraphics[width=1.0\linewidth]{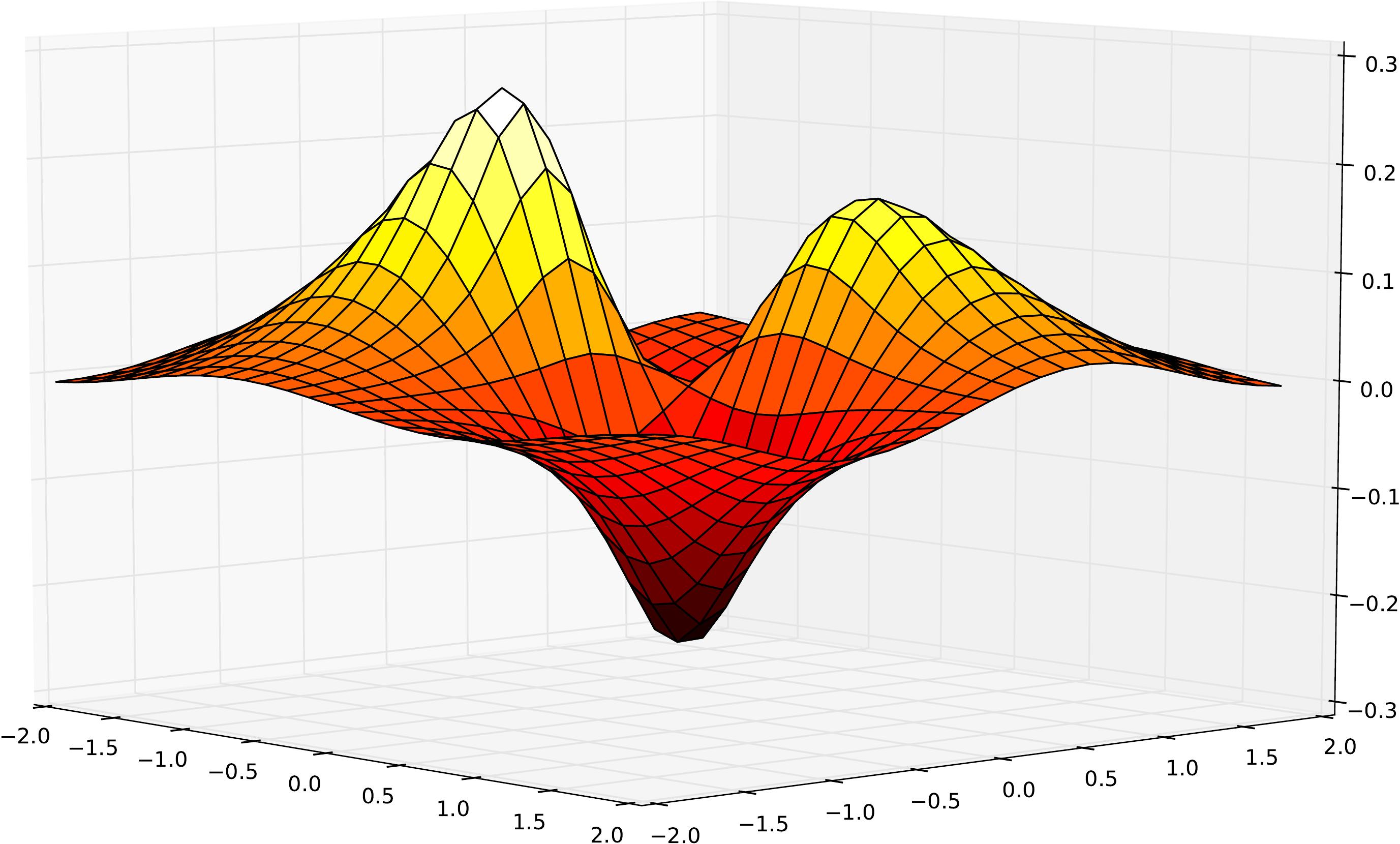}
  \caption{Nonlinear optimization}
  \label{fig:sub2}
\end{subfigure}
\caption{Optimization problems of different degree of complexity are described above. Computational  decisions  that depend on complex optimization algorithms take longer time to compute, affecting the mission effectiveness.}
\label{fig:optimization}
\end{figure}

Given each mission has a computational requirement that must be satisfied, a lot of effort is expended to find computationally efficient algorithms. One gauge of efficiency is based on the computational complexity  of an algorithm which is generally expressed in terms of time complexity that describes how long it takes for an algorithm to compute an answer using limited number of resources. The time complexity will affect the efficiency of  applications using those algorithms. In mission-oriented tactical environments, computational efficiency  in terms of flops/watt and computational speeds to match mission requirements  is a very important factor in determining the fitness of an application for mission deployment. The requirements include computational efficiency of each platform in regards to available resources and constraints. Network specific constraints also need to be taken into consideration when assessing the efficiency of distributed computation. Furthermore, the complexity of input tasks and the computational decision-making requirements will increase, having variable computational cost on each platform, and often these functions need to be optimized\cite{fazel2005network, lee2005non, nygren2010akamai}. The type of optimization algorithms  will also have an effect on the time to compute. For instance, optimizations can be linear or nonlinear, and nonlinear optimization can be further classified into concave and convex optimizations based on the complexity (see Fig.~\ref{fig:optimization}). Linear optimization functions run faster than nonlinear optimization functions as linear programming isolates computation to just the vertices of the available parameter space while nonlinear optimization must continuously explore the parameter space. In addition to the complexity of the computational tasks, resource constraints will severely affect the  ability of tactical computing platforms to complete assigned tasks in the desired time.  Some of the mission critical applications demand well-defined execution times. Mission optimized computations combined with automated intelligence have been shown to be efficient and will achieve the desired results in mission time\cite{spang2017mon, capacity-bounds-multipath-networks, Afergan:2005:EPB:1251522.1251523}. The mission requirements define the need for optimized algorithms, which are a subset of  the deterministic polynomial time complexity class $P$ of problems, where we define it as the $M$ class for mission ready algorithms. 

Context-aware adaptive  computational frameworks are remarkably flexible to support the class $M\subseteq P$  of languages using heterogeneous underlying hardware implementation details from multiple applications using non-intrusive methods under resource constrained tactical environments. Mission optimized adaptive computational framework  will improve the execution efficiency of the mission applications. Adaptive computing framework  will be highly useful in deploying intelligent  tactical computing platforms in smart cities that present complex environments with constrained resources and heterogeneity with optimal efficiency. It will also  enable the computing platforms to adapt to the communication related constraints like bandwidth and  network reachability. In tactical environments, adaptive computing must make decisions  regarding local or remote computing for solving a computational problem presented to it.  Computational offloading is influenced by not only the computational complexity of the problem, but it also depends upon the network resource constraints like the link quality and communication cost for the offloading.

With military missions,  efficiency of computing platforms  are affected by  resource constraints  in the field. Availability of the limited amount of energy and battery power, the timing of communication lines for crucial decision-making, and the efficiency of available central processing power (CPU) that are accessible in real time are immensely crucial aspects in order to carry out a successful computation required by the mission. By focusing on military operations, we will analyze to optimize necessary computations in a mission where computations are computable in polynomial time and investigate when such computations can be optimized, reduced in the usage of CPU power, and successfully executed within mission time. Such study will require us to develop a heterogeneous (non)linear platform to reduce current state-of-the-art computational complexity to the complexity sufficient for military applications. 

In this manuscript, we study a new class of polynomial time computational complexity called  $M$ to satisfy the mission requirements in order to understand the computability of a given algorithm in mission time. Mission times can vary based on the mission objectives but deployability of a given algorithm is tied to its ability to complete its computation in that particular mission time. In the section on computational-complexity, we give a brief background on computational complexity, while in the section on adaptive mission computation, we define this new $M$ class of algorithms and we specifically look at how algorithms can be optimized to be placed from the $P$ class of algorithms into the $M$ class by using an adaptive computing framework. 
Any adaptive computing framework must consider the constraints which we outline in the section regarding the effect of resources on computational efficiency. Next we apply those constraints towards a constraints-aware distributed computing framework, and then we give an example where computational jobs assigned to a cluster of local machines may fail and how the distributed computing platform reacts to these failures. In the final section, we summarize our construction of the mission class of problems and how focusing on this subset will enhance the performance of tactical mobile computing platforms. 

\section{Computational Complexity}\label{section:computational-complexity}
Computational complexity has a direct impact on mission ready algorithms. For example, RSA encryption is simple if given the public key, but RSA decryption is difficult without knowledge of the private key. As such, asymmetric encryption can be performed easily in the field, but mounting a brute force attack on an adversary using RSA encryption is not even considered in real-time applications. This variation in complexity is why we group problems into two types: class $P$, and class $NP$\cite{bossaerts2017computational}.

We denote $P$ as the class of questions for which some algorithm can provide an answer, and thus solve the language, in polynomial time. So in the mathematics literature, we often refer to $P$ as {\em deterministic polynomial time complexity class}. The class $P$ of problems are decidable in polynomial time on a deterministic single-tape Turing machine: $P=\bigcup_{k\geq 0}\text{Time}(n^k)$, 
and such problems are simple for computers to solve, all within a reasonable amount of time. One example includes determining whether or not a word $w$ is a member of the language $L=\{0^k1^k :k\geq 0\}$. 
Another example includes the problem to determine whether a directed path exists from vertex $s$ to vertex $t$ in a directed graph, i.e., given a directed graph $G$, define 
\begin{align}
\begin{split}
&\text{PATH}(G,s,t) :=\\
& \{ \langle G,s,t\rangle : G \text{ has a directed path from }s \mbox{ to }t \}. 
\end{split}
\end{align}
This problem is known as {\em directed $s$-$t$ connectivity}, and it is a classical result that $\text{PATH}$ is indeed in the class $P$\cite{barnes1998sublinear}. 
In fact, reasonable deterministic computational models are polynomially equivalent, i.e., one such model can simulate another model with only a polynomial increase in running time. 

We denote $NP$ as the complexity class for which an algorithm can provide a solution in polynomial time with a non-deterministic Turing machine. We thus refer to $NP$ as {\em nondeterministic polynomial time complexity class}. Such questions are ones with solutions that can be {\em verified} in polynomial time using a deterministic Turing machine. By definition, $NP$ is the class of languages that are decidable in polynomial time on a nondeterministic Turing machine, i.e., $NP=\bigcup_{k\geq 0}\text{NTime}(n^k)$, where 
\begin{align}
\begin{split}
&\text{NTime}(t(n)) := \{ L : L\mbox{ is a language decided by an}\\
& \mathcal{O}(t(n))\mbox{-time nondeterministic Turing machine}\}. 
\end{split}
\end{align}

In fact, many questions in the class $P$ can be changed ever so slightly to then be placed in the class $NP$. For example, the previous directed path question that connects two vertices can be placed into the $NP$ class by instead considering the question of whether or not a Hamiltonian path connects two vertices, where a {\em Hamiltonian path} in a directed graph $G$ is a directed path that passes through each vertex exactly once. That is, given a directed graph $G$, we define a {\em Hamiltonian path} as
\begin{align}
\begin{split}
&\text{HamPATH}(G,s,t) = \{ \langle G,s,t\rangle: G \mbox{ is a directed}\\
& \mbox{graph with a Hamiltonian path from }s \mbox{ to }t \}. 
\end{split}
\end{align}

One can easily obtain an {\em exponential} time algorithm for the $\text{HamPATH}$  problem by a brute-force approach that checks all possible permutations of vertices (if there are $n$ vertices, then there are $n!$ permutations to check, and we only need to verify that a potential path is Hamiltonian). Although it is well-known that $\text{HamPATH}$ is in $NP$,\cite{gurevich1987expected, itai1982hamilton} 
it remains an open problem on determining whether or not $\text{HamPATH}$ is actually solvable in polynomial time. 
We thus have the containment $P\subseteq NP$ of classes of languages. 

\begin{figure}
\centerline{\includegraphics[scale=0.88]{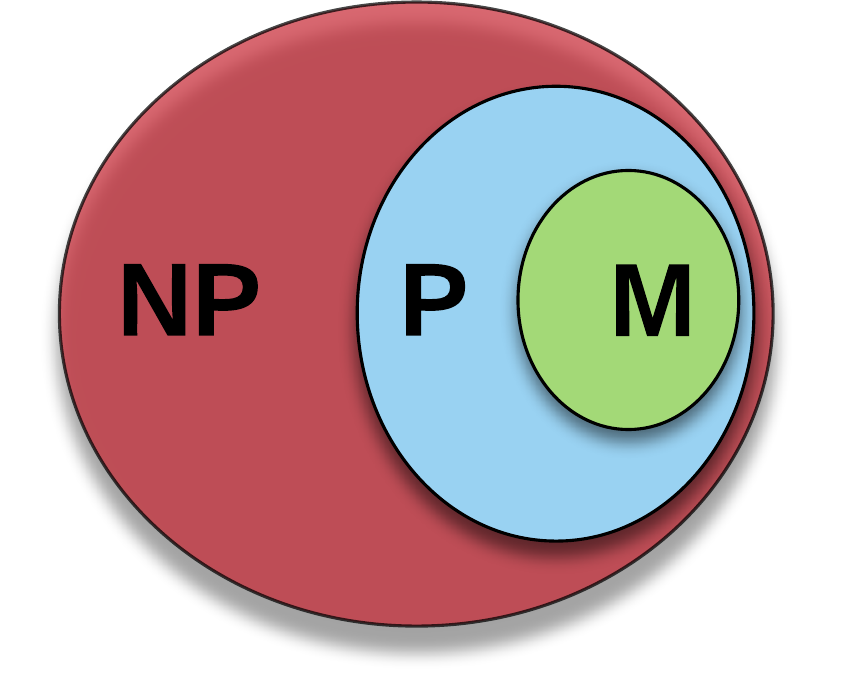}}
\caption{We say $P$ (deterministic polynomial time complexity class) is the class of languages for which some algorithm can solve the question in polynomial time while 
$NP$ (nondeterministic polynomial time complexity class) is the class of languages for which an algorithm may be very difficult to find, but if provided an answer, then it can be {\em verified} in polynomial time. The subclass $M$ consists of functions in $P$ for which adequate  computational resources  may not be available in the tactical environment to complete mission computation but for which a significant proportion of the language can be completed within the confinement of limited resources, such as technology and time, so that the validity of a solution for the language may correctly be deduced decisively.}
\label{fig:mission_comp} 
\end{figure}

In the study of computability, $P=NP$ problem is one of the very important and interesting open problems in mathematics and in theoretical computer science with deep ramifications in cryptography\cite{baker1975relativizations, marks2016universality}, algorithms\cite{cook2006p, yan2002number, shor2004progress}, artificial intelligence\cite{cook2006p, buss2012towards}, game theory\cite{sahni1974computationally}, economics\cite{vives1984duopoly}, to name a few. Intuitively, the $P=NP$ conjecture is an investigation of whether every problem whose solution can be quickly verified in polynomial time can also be solved fairly quickly in polynomial time. So it remains an open problem to show whether or not the following containment $NP\stackrel{?}{\subseteq} P$ of languages holds. 
  
  There is however a class of (optimization) problems that lie in the class $NP$, so a polynomial time reduction to their complexity may be difficult to construct with current tools.
  
Now, suppose a full high performance computing (HPC) may not be available on the mission field but for which a significant proportion of the language being computed within a confined time frame is sufficient to determine the validity of the language, and thus critical decisions may need to be immediately determined on site. Because of the importance and a critical significance of such problems, we define a subclass of problems in $P$ in the section on adaptive mission computation.

\section{Adaptive Mission Computation}
\label{section:related-work}
With adaptive computing in mind~\cite{pozueco2013adaptable, urbieta2017adaptive}, the class $M$ of functions can be described as having the ability to manage time constraint and the ability to complete a computation in the class $P$ as we minimize the usage of limited HPC resources, 
such as energy, time, and memory, 
while simultaneously attempt to maximize computational power (CPU) and efficiency to effectively work around computational constraints for a heterogeneous computational platform for the class $M$ of languages. 
The timing of such computations may be optimized by using a convex or concave  (non)linear platform, and by investigating how computations perform in a mission, the limits of computational complexity and computational capacity may be precisely described.

We are interested in a mission-focused problem $\mathcal{P}$, which is a relation from a set $I$ of instances (input) to a set $S$ of solutions, where deterministic and approximation algorithms exist, i.e., 
$\mathcal{P}\subseteq I\times S$, $\mathcal{P}\in P$, and $\mathcal{P}$ has an algorithm that consistently returns a feasible, approximated solution, which is characterized by its distance from its value to the optimal solution. 
Thus the class $M$ is defined as the following: 
\begin{definition}\label{defn:mission-computable-class}
Given the polynomial time complexity class  
\begin{equation}
P =\{ \mathcal{P}: \forall x\in I\: \exists \: y \in S \ni (x,y)\in \mathcal{P}\},  
\end{equation}
we define 
{\em mission computable polynomial time complexity class}  
as 
\begin{equation}
M =\{ \mathcal{P}\in P: \forall x\in I \: \exists \mbox{ approx. soln. }
\widetilde{y} \ni (x,\widetilde{y})\in \mathcal{P}\}. 
\end{equation}
\end{definition} 

Note that the approximated solution $\widetilde{y}$ does not need to be in the solution set $S$, but given any $\varepsilon > 0$, $\widetilde{y}$ must satisfy $d(y,\widetilde{y})< \varepsilon$, where $d$ is an intrinsic notion of a distance between the two solutions $y$ and $\widetilde{y}$, which depends on the problem being considered. 

Our class $M$ of computational complexity is aimed at solving computational problems  in  mission time and it is restricted to adaptive computing framework  (see Fig.~\ref{fig:mission_comp}).

% we have added more quantifiers 

Before we further discuss approximated solutions, we will give some basic conditions about $\mathcal{P}$.  In order to solve the language $\mathcal{P}$, one needs to recognize if $(x,y)\in \mathcal{P}$. One then needs to construct (using a deterministic algorithm with polynomial time complexity) that for each $x$, there exists $y$ such that $(x,y)\in \mathcal{P}$. Finally, we need to optimize the problem $\mathcal{P}$ in such a way that for each instance $x$, find the best (and efficient) solution $\widetilde{y}$ such that $(x,\widetilde{y})\in \mathcal{P}$.

A feasible solution is an approximate solution, and such solutions are classified by the value of its distance from the ideal (optimal) solution. Thus, the ratio of the rough solution to the optimal solution is determined by the input size, the growth of the performance algorithm, and time limitations. Some current approximation algorithms include linear programming, dynamic programming, local search, randomized algorithm, and heuristic algorithm. Linear programming formulates the language as a linear model, dynamic programming constructs a solution from optimal solutions to sub-problems of the original language, a local search algorithm looks for a better neighboring solution when given a solution (this algorithm {\em deforms} an input until an improved and preferred optimal solution is found), randomized algorithm embeds and executes a random decision generator,  and heuristic algorithm is encoded by exploratory learning strategies that offer no guarantee of a more suitable solution. 

The mission class $M$ is the subclass of optimization problems that are solvable by a polynomial time and approximate algorithm, in a finite sequence of steps (whose order is bounded by pre-determined complexity), that computes a comparable result when given an instance from the set of inputs. The algorithm cost for the number of operations (time complexity) and the storage space (space complexity) are in the order of polynomial time $\text{Time}(n^k)$ for some $k$. 

\section{Effect of Resources on Computational Efficiency}
\label{section:resource-effects}
For smooth operations, adaptive computations require allocation algorithms and technology congestion protocols using behavior models along with network efficiency fairness characterization. These algorithms need to make an optimal decision which can be characterized as an optimization problem that is either linear or nonlinear. As such, objective functions together with a set of constraint inequalities are often used to broaden the scope for a successful military operation. Computational constraints like power, memory, size, storage and CPU  influence the performance of a computational platform in contested and congested environments, and are driving the need for a constraints-aware adaptive computation framework. Such a framework will change the computational behaviors of the platforms in response to available resources and the complexity of the input computational tasks. For example, a computational problem can be solved in a distributed manner in order to optimize available computational resources among different computational platforms. However, when a local computing platform is incapable of performing a required computation due to lack of resources or due to lack of required software, the computations are offloaded to a remote computer capable of computing the problem and provide the solution. The network related constraints like signal strength, bandwidth and energy required for offloading are important factors in determining if a remote computation makes sense. Any adaptive computing framework will need to take these factors into consideration when making decisions.

\section{Constraints-aware distributed computing}
\label{section:computing-with-constraints}
One foundation for adaptive computing is a constraints-aware distributed computing algorithm. Given too many jobs assigned to an array of cores, the algorithm can be programmed to minimize the number of failed jobs; such an algorithm is inherently designed to allow for failure, but it is this failure that allows for feedback to the calling applications. For example, a security camera might be running at $60$ Hz, but the image analysis at $60$ Hz would consume more than the local resources. The application dedicated to analyzing the images would give a time-to-complete restriction on each frame, but a local optimizer would determine that most of the frames will fail the time restrictions and terminate most jobs immediately while informing the calling application about the dropped jobs.

The logic for a single machine with a single core can be programmed using integer programming. Given we want to maximize the total computations done on a single machine, our decision variables can be an ordered list that is represented as a binary matrix $b_{i,j}$ (see Eq.~\ref{eq:decisionVariables}), where $i$ represents the $i$-th job that will be computed in order on the host machine, while $j$ represents whether or not the $j$-th image is assigned to the $i$-th job. Given the decision variables, the objective function can be used to maximize the computations done on the local machine by doing a weighted sum (see the first term of Eq.~\ref{eq:objectiveFunction}). Constraints are added to prevent multiple images being assigned to the same job as in Eq.~\ref{eq:constraint1} or to prevent the same image being applied to multiple jobs as in Eq.~\ref{eq:constraint1b}, while the time constraints in Eq.~\ref{eq:constraint2} are put in place to prevent images from being scheduled too far into the future on the local machine. Integer programming must work within the constraints. 
However, we expect some jobs to fail to be assigned, so a null image with zero computation time and zero constraints is assigned to the $0$-th image. Furthermore, a compact list with all null jobs appearing at the end is preferable to limit degeneracy, so the objective function has an additional term for the zeroth null image appearing later in the list (see the second term of Eq.~\ref{eq:objectiveFunction}). Feeding the variables, objective function, and constraints into an integer optimizer will show which jobs should be abandoned on the local machine. These abandoned jobs can then be scheduled on a distributed computing platform.

%\begin{equation}\label{eq:decisionVariables}
\begin{align}\label{eq:decisionVariables}
\begin{split}
i&=\text{job index},  \hspace{4mm}
j=\text{image index},\\
b_{i,j}&= 
	\begin{cases} 
0 &\text{if $j$-th image not assigned to job $i$}, \\
1 & \text{if $j$-th image assigned to job $i$}, 
	\end{cases}
\end{split}
\end{align}
%\end{equation}

\begin{align}\label{eq:helperFunctions}
\begin{split}
T_j&=\text{analysis time for the $j$-th image}, \\
R_j&=\text{time the $j$-th image must be processed}, 
\end{split}
\end{align}

objective function:
\begin{equation}\label{eq:objectiveFunction}
\sum_{i,j\neq0}b_{i,j}\cdot T_j \text{ }+10^{-5}\cdot\sum_i i\cdot b_{i,0}, 
\end{equation}

constraint: one image per job
\begin{equation}\label{eq:constraint1}
\forall i \hspace{4mm} \text{ } \sum_{j} b_{i,j}=1,   
\end{equation}

constraint: each image is processed at most once
\begin{equation}\label{eq:constraint1b}
\forall j\neq 0  \hspace{4mm} \text{ } \sum_{i} b_{i,j}\leq 1,   
\end{equation}

constraint: image processed on time:
\begin{equation}\label{eq:constraint2}
\forall i \text{ } \hspace{4mm} \sum_{k\leq i,j} b_{k,j}\cdot T_j\leq \sum_j b_{i,j}\cdot R_{j}.  
\end{equation}

Once the abandoned jobs are handed to a distributed computing platform, a global optimizer can determine which jobs should be scheduled or rescheduled. The distributed computing platform is really just a bunch of local machines connected together with each machine competing for the same resources, so again each application must be willing to accept that some jobs will be dropped. This means each machine must throttle its own computational usage through tagging its own jobs with the appropriate priority level. In the scenario of image analysis, every tenth frame can be tagged as a high priority job relative to the calling machine. Although the local optimizer outlined previously is fairly trivial, the same framework can be expanded to include a number of cores ($c$) and a number of machines ($m$). The decision variables are again just an expansion of the binary matrix, from $b_{i,j} $ to $b_{m,c,i,j}$. Extra constraints can be added to satisfy transfer times, power used, or even the amount of random-access memory (RAM) used. The objective function can be tweaked to emphasize particular features such as computation time or energyConsumed/FLOP.

The common theme between the local optimizer and the global optimizer is that the computational size must be varied to fit the computational resources. This is just a foundation for adaptive computing as the next step is to vary the allocation of resources to competing applications. For instance image analysis across multiple security cameras can be considered a single image analysis application, while another application could be attempting to apply machine learning to detect hostile agents from the data provided by the security cameras. The overall goal is to maximize mission effectiveness by varying the resources to each application as the needs change. For instance, image analysis in hostile zones should be given priority access to local resources, while the machine learning project should be given loose timing constraints to allow the jobs to be scheduled on HPC machines located far away from the hostile zone.

\section{Discussion}
\label{section:discussion}

Tactical computing platforms are mostly mobile with limited computing and communication resources. Data processing and problem solving tasks are time sensitive and  their speed depends upon the available computational resources. In order to accomplish  mission computation goals, the platforms and the algorithms must be optimized to the mission requirements. 

We have described a new class of computational complexity class $M$  that is a subset of the polynomial time class $P$ in order to address a class of problems that needs to be computed in mission time. Mission times are determined by context in which the computation is used and the completion time of that task to determine the usefulness of the computation. As described in the section on adaptive mission computation, all the polynomial class of computational tasks that can be completed in mission time will fall into the computational complexity class $M$, and they are mission ready. Polynomial problems that cannot be computed in mission time will require additional optimization until they satisfy mission requirements.

Defining a new computational complexity class will allow us to define the computational requirements for any applications and algorithms to be mission ready. Algorithms that  can process input data in polynomial time might still fail to complete  when deployed in tactical environments due to resource constraints. Reducing such a $P$ class of algorithms to $M$ class through optimization and heuristics  will enable us to optimize a given set for $P$ class of problems to tactical environments. 

In our future work, we will test a variety of mission deployable applications  to determine their  mission readiness without additional optimization and adjustments to their code.
  
\section{Acknowledgements}

This work is supported by a research collaboration between the U.S. Army Research Laboratory and U.S. Military Academy. 
M.S.I. is partially supported by National Research Council Research Associateship Programs.

%\begin{verbatim}
%%Harvard (name/date)
%\bibliographystyle{SageH}
%%Vancouver (numbered)
%\bibliographystyle{Sagev}
\bibliographystyle{IEEEtran}
\bibliography{im-dasari}

\begin{thebibliography}{10}
\providecommand{\url}[1]{#1}
\csname url@samestyle\endcsname
\providecommand{\newblock}{\relax}
\providecommand{\bibinfo}[2]{#2}
\providecommand{\BIBentrySTDinterwordspacing}{\spaceskip=0pt\relax}
\providecommand{\BIBentryALTinterwordstretchfactor}{4}
\providecommand{\BIBentryALTinterwordspacing}{\spaceskip=\fontdimen2\font plus
\BIBentryALTinterwordstretchfactor\fontdimen3\font minus
  \fontdimen4\font\relax}
\providecommand{\BIBforeignlanguage}[2]{{%
\expandafter\ifx\csname l@#1\endcsname\relax
\typeout{** WARNING: IEEEtran.bst: No hyphenation pattern has been}%
\typeout{** loaded for the language `#1'. Using the pattern for}%
\typeout{** the default language instead.}%
\else
\language=\csname l@#1\endcsname
\fi
#2}}
\providecommand{\BIBdecl}{\relax}
\BIBdecl

\bibitem{fazel2005network}
M.~Fazel and M.~Chiang, ``Network utility maximization with nonconcave
  utilities using sum-of-squares method,'' in \emph{Decision and Control, 2005
  and 2005 European Control Conference. CDC-ECC'05. 44th IEEE Conference
  on}.\hskip 1em plus 0.5em minus 0.4em\relax IEEE, 2005, pp. 1867--1874.

\bibitem{lee2005non}
J.-W. Lee, R.~R. Mazumdar, and N.~B. Shroff, ``Non-convex optimization and rate
  control for multi-class services in the internet,'' \emph{IEEE/ACM
  transactions on networking}, vol.~13, no.~4, pp. 827--840, 2005.

\bibitem{nygren2010akamai}
E.~Nygren, R.~K. Sitaraman, and J.~Sun, ``The {A}kamai network: a platform for
  high-performance internet applications,'' \emph{ACM SIGOPS Operating Systems
  Review}, vol.~44, no.~3, pp. 2--19, 2010.

\bibitem{spang2017mon}
B.~Spang, A.~Sabnis, R.~Sitaraman, D.~Towsley, and B.~DeCleene, ``{MON}:
  {M}ission-optimized overlay networks,'' in \emph{INFOCOM 2017-IEEE Conference
  on Computer Communications, IEEE}.\hskip 1em plus 0.5em minus 0.4em\relax
  IEEE, 2017, pp. 1--9.

\bibitem{capacity-bounds-multipath-networks}
\BIBentryALTinterwordspacing
A.~Bejan, R.~Gibbens, R.~Hancock, and D.~Towsley, ``Capacity bounds and
  robustness in multipath networks,'' in \emph{Proceedings of the 9th EAI
  International Conference on Performance Evaluation Methodologies and Tools},
  ser. VALUETOOLS'15.\hskip 1em plus 0.5em minus 0.4em\relax ICST, Brussels,
  Belgium, Belgium: ICST (Institute for Computer Sciences, Social-Informatics
  and Telecommunications Engineering), 2016, pp. 9--16. [Online]. Available:
  \url{http://dx.doi.org/10.4108/eai.14-12-2015.2262625}
\BIBentrySTDinterwordspacing

\bibitem{Afergan:2005:EPB:1251522.1251523}
\BIBentryALTinterwordspacing
M.~Afergan, J.~Wein, and A.~LaMeyer, ``Experience with some principles for
  building an internet-scale reliable system,'' in \emph{Proceedings of the 2Nd
  Conference on Real, Large Distributed Systems - Volume 2}, ser.
  WORLDS'05.\hskip 1em plus 0.5em minus 0.4em\relax Berkeley, CA, USA: USENIX
  Association, 2005, pp. 1--6. [Online]. Available:
  \url{http://dl.acm.org/citation.cfm?id=1251522.1251523}
\BIBentrySTDinterwordspacing

\bibitem{bossaerts2017computational}
P.~Bossaerts and C.~Murawski, ``Computational complexity and human
  decision-making,'' \emph{Trends in cognitive sciences}, vol.~21, no.~12, pp.
  917--929, 2017.

\bibitem{barnes1998sublinear}
G.~Barnes, J.~F. Buss, W.~L. Ruzzo, and B.~Schieber, ``A sublinear space,
  polynomial time algorithm for directed $s$-$t$ connectivity,'' \emph{SIAM
  Journal on Computing}, vol.~27, no.~5, pp. 1273--1282, 1998.

\bibitem{gurevich1987expected}
Y.~Gurevich and S.~Shelah, ``Expected computation time for {H}amiltonian path
  problem,'' \emph{SIAM Journal on Computing}, vol.~16, no.~3, pp. 486--502,
  1987.

\bibitem{itai1982hamilton}
A.~Itai, C.~H. Papadimitriou, and J.~L. Szwarcfiter, ``Hamilton paths in grid
  graphs,'' \emph{SIAM Journal on Computing}, vol.~11, no.~4, pp. 676--686,
  1982.

\bibitem{baker1975relativizations}
T.~Baker, J.~Gill, and R.~Solovay, ``Relativizations of the {P}=?{NP}
  question,'' \emph{SIAM Journal on computing}, vol.~4, no.~4, pp. 431--442,
  1975.

\bibitem{marks2016universality}
A.~Marks, ``The universality of polynomial time {T}uring equivalence,''
  \emph{Mathematical Structures in Computer Science}, pp. 1--9, 2016.

\bibitem{cook2006p}
S.~Cook, ``The {P} versus {NP} problem,'' \emph{The millennium prize problems},
  pp. 87--104, 2006.

\bibitem{yan2002number}
S.~Y. Yan, \emph{Number theory for computing}.\hskip 1em plus 0.5em minus
  0.4em\relax Springer Science \& Business Media, 2002.

\bibitem{shor2004progress}
P.~W. Shor, ``Progress in quantum algorithms,'' \emph{Quantum information
  processing}, vol.~3, no. 1-5, pp. 5--13, 2004.

\bibitem{buss2012towards}
S.~R. Buss, ``Towards {NP--P} via proof complexity and search,'' \emph{Annals
  of Pure and Applied Logic}, vol. 163, no.~7, pp. 906--917, 2012.

\bibitem{sahni1974computationally}
S.~Sahni, ``Computationally related problems,'' \emph{SIAM Journal on
  Computing}, vol.~3, no.~4, pp. 262--279, 1974.

\bibitem{vives1984duopoly}
X.~Vives, ``Duopoly information equilibrium: {C}ournot and {B}ertrand,''
  \emph{Journal of economic theory}, vol.~34, no.~1, pp. 71--94, 1984.

\bibitem{pozueco2013adaptable}
L.~Pozueco, X.~G. Pa{\~n}eda, R.~Garc{\'\i}a, D.~Melendi, and S.~Cabrero,
  ``Adaptable system based on scalable video coding for high-quality video
  service,'' \emph{Computers \& Electrical Engineering}, vol.~39, no.~3, pp.
  775--789, 2013.

\bibitem{urbieta2017adaptive}
A.~Urbieta, A.~Gonz{\'a}lez-Beltr{\'a}n, S.~B. Mokhtar, M.~A. Hossain, and
  L.~Capra, ``Adaptive and context-aware service composition for {I}o{T}-based
  smart cities,'' \emph{Future Generation Computer Systems}, vol.~76, pp.
  262--274, 2017.

\end{thebibliography}
%\end{verbatim}

%\section{Support for \textsf{\journalclass}}
%We offer on-line support to participating authors. Please contact
%us via e-mail at \dots
%
%We would welcome any feedback, positive or otherwise, on your
%experiences of using \textsf{\journalclass}.

%\section{Copyright statement}
%Please  be  aware that the use of  this \LaTeXe\ class file is
%governed by the following conditions.
%
%\subsection{Copyright}
%Copyright \copyright\ \volumeyear\ SAGE Publications Ltd,
%1 Oliver's Yard, 55 City Road, London, EC1Y~1SP, UK. All
%rights reserved.
%
%\subsection{Rules of use}
%This class file is made available for use by authors who wish to
%prepare an article for publication in a \textit{SAGE Publications} journal.
%The user may not exploit any
%part of the class file commercially.
%
%This class file is provided on an \textit{as is}  basis, without
%warranties of any kind, either express or implied, including but
%not limited to warranties of title, or implied  warranties of
%merchantablility or fitness for a particular purpose. There will
%be no duty on the author[s] of the software or SAGE Publications Ltd
%to correct any errors or defects in the software. Any
%statutory  rights you may have remain unaffected by your
%acceptance of these rules of use.
%
%\begin{acks}
%This class file was developed by Sunrise Setting Ltd,
%Brixham, Devon, UK.\\
%Website: \url{http://www.sunrise-setting.co.uk}
%\end{acks}

\end{document}